\documentclass[a4paper,leqno]{amsart}

\usepackage{amsmath,amssymb,stmaryrd}
\usepackage{amsthm}

\newcommand{\R}{\mathbb{R}}
\newcommand{\N}{\mathbb{N}}
\newcommand{\Z}{\mathbb{Z}}
\renewcommand{\d}{\partial}

\newcommand{\pare}[1]{\left(#1\right)}
\newcommand{\quadre}[1]{\left[#1\right]}
\newcommand{\curly}[1]{\left\{#1\right\}}

\newcommand{\modu}[1]{\left|#1\right|}
\newcommand{\eps}{\varepsilon}
\newcommand{\cfun}{\mathcal C}

\newcommand{\dist}{\mathcal D}

\newtheorem{theorem}{Theorem}[section]
\newtheorem{lemma}[theorem]{Lemma}

\newtheorem{definition}[theorem]{Definition}

\newtheorem{corollary}[theorem]{Corollary}
\newtheorem{remark}[theorem]{Remark}

\numberwithin{equation}{section}

\title[On the XFEL Schr\"odinger Equation]{On the XFEL Schr\"odinger Equation: Highly Oscillatory Magnetic Potentials and Time Averaging}
\author[P. Antonelli]{Paolo Antonelli}
\address[P. Antonelli]{CEREMADE, UniversitŽ de Paris-Dauphine\\
Place du MarŽchal De Lattre De Tassigny\\
75775 PARIS CEDEX 16 - FRANCE}
\email{antonelli.p@gmail.com}
\author[A. Athanassoulis]{Agisillaos Athanassoulis}
\address[A. Athanassoulis]{Department of Applied Mathematics\\
University of Crete, Heraklion 71409\\Greece}
\email{athanassoulis@tem.uoc.gr}
\author[H. Hajaiej]{Hichem Hajaiej}
\address[H. Hajaiej]{King Saud University (KSU),
College of Science\\
Riyadh,
Kingdom of Saudi Arabia}
\email{hichem.hajaiej@gmail.com}
\author[P. Markowich]{Peter Markowich}
\address[P. Markowich]{King Abdullah University of Science and Technology (KAUST)\\
MCSE Division\\
Thuwal 23955-6900, Kingdom of Saudi Arabia}
\email{Peter.Markowich@kaust.edu.sa}

\date{\today}

\begin{document}
\maketitle
\begin{abstract}
We analyse a nonlinear Schr\"odinger equation for the time-evolution of the wave function of an electron beam, interacting selfconsistently through a Hartree-Fock nonlinearity and through the repulsive Coulomb interaction of an atomic nucleus.
The electrons are supposed to move under the action of a time dependent, rapidly periodically oscillating electromagnetic potential. This can be considered a simplified effective single particle model for an X-ray Free Electron Laser (XFEL). We prove the existence and uniqueness for the Cauchy problem and the convergence of wave-functions to corresponding solutions of a Schr\"odinger equation with a time-averaged Coulomb potential in the high frequency limit for the oscillations of the electromagnetic potential.
\end{abstract}

\section{Introduction and Statement of the Main Results}\label{sect:first}
In this paper we investigate the following nonlinear Schr\"odinger equation
\begin{equation}\label{eq:xfel}
i\hbar\d_t\psi=(i\hbar\nabla-A)^2\psi+c\frac{1}{|x|}\psi+C_1(|\cdot|^{-1}\ast|\psi|^2)\psi
-a|\psi|^\sigma\psi.
\end{equation}
The coefficients $c, C_1, a$ and the exponent $\sigma$ are assumed to be nonnegative and $\hbar$ is supposed to be a scaled version of the Planck constant, which - w.l.o.g. - shall be set equal to 1 in the sequel. A solution $\psi$ of this Schr\"odinger equation can be considered as the wavefunction of an electron beam, interacting self-consistently through the repulsive Coulomb (Hartree) force with strength $C_1$, the attractive local Fock approximation with strength $a$ (later on we shall comment on the exponent $\sigma$) and interacting repulsively with an atomic nucleus, located at the origin, of interaction strength $c$. The vectorfield $A$ represents an external electromagnetic field, which we shall assume to depend on time $t$ only (not on position $x$). Clearly, this Schr\"odinger equation is time-reversible, but for the sake of notational simplicity we consider $t>0$. For physical reasons we shall only consider the three-dimensional case here, i.e. the spatial variable $x$ is assumed to be in $\R^3$. 
\newline
Nevertheless, because of the sole dependence of $A$ on time, by a simple change of coordinates and a phase shift, we see that equation \eqref{eq:xfel} can be transformed into a similar nonlinear Schr\"odinger equation, where the electromagnetic Laplacian is replaced by the standard one, but on the other hand a time-dependent Coulomb potential appears.
Indeed, by defining
\begin{equation}\label{eq:psi_to_u}
u(t, x)=\psi(t, x+b(t))e^{i\int_0^t|A(s)|^2ds},
\end{equation}
where $b(t)=2\int_0^tA(s)ds$, then we can see that $u$ satisfies 
\begin{equation}\label{eq:xfel_2}
i\d_tu=-\Delta u+Vu+C_1(|\cdot|^{-1}\ast|u|^2)u-a|u|^\sigma u,
\end{equation}
where now the potential is given by
\begin{equation*}
V(t, x)=\frac{c}{|x-b(t)|}.
\end{equation*}
In this paper we are interested in studying the case when $A(t)$ is rapidly oscillating and we investigate the asymptotic behaviour of solutions of \eqref{eq:xfel_2} in the highly oscillating regime. Then the equation \eqref{eq:xfel} can be considered as a model for XFEL (X-Ray Free Electron Laser), cf. \cite{FR}.
\newline
As an example we can think of $b(t)=\vec e\sin(\omega t)$, where $\omega\gg1$ is the oscillation frequency, and $\vec e$ is a constant vector in $\R^3$, but as we will show this can be extended to the case where the field $b$ can be written as
\begin{equation}\label{eq:125}
b(t)=\vec e(t)f(\omega t),
\end{equation}
where $\vec e:\R\to\R^3$ is a smooth vector field and $f$ is an arbitrary continuous, $2\pi-$periodic function.
\newline
To this end, we shall point out the $\omega$ dependence of functions with a superscript, 
$b^\omega(t)=\vec e(t)f(\omega t)$, and
\begin{equation}\label{eq:V_omega}
V^\omega(t, x)=\frac{c}{|x-b^\omega(t)|}.
\end{equation}
We will then study solutions of the Cauchy problem
\begin{equation}\label{eq:nls_omega_intro}
\left\{\begin{array}{l}
i\d_tu^\omega=-\Delta u^\omega+V^\omega u^\omega+C_1(|\cdot|^{-1}\ast|u^\omega|^2)u^\omega-a|u^\omega|^\sigma u^\omega\\
u^\omega(0)=u_0,
\end{array}\right.
\end{equation}
and their convergence to solutions of the \emph{averaged} equation
\begin{equation}\label{eq:nls_aver_intro}
\left\{\begin{array}{l}
i\d_tu=-\Delta u+\langle V\rangle u+C_1(|\cdot|^{-1}\ast|u|^2)u-a|u|^\sigma u\\
u(0)=u_0,
\end{array}\right.
\end{equation}
where $\langle V\rangle$ is the limiting potential and is given (see Section \ref{sect:intro} for details) by
\begin{equation}\label{eq:V_aver}
\langle V\rangle(t, x):=\int_0^1\frac{c}{|x-\vec e(t)\sin(2\pi\omega\tau)|}d\tau.
\end{equation}
The main theorem we will prove in this paper is the following one
\begin{theorem}\label{thm:main}
Let $0<\sigma<4/3$, $u_0\in L^2(\R^3)$, $u^\omega, u\in\cfun(\R; L^2(\R^3))$ be the unique global solutions of 
\eqref{eq:nls_omega_intro}, \eqref{eq:nls_aver_intro}, respectively (see Theorem \ref{thm:mass_subcrit} below). Then 
for each finite time $0<T<\infty$ and for each admissible Strichartz index pair $(q, r)$, we have
\begin{equation*}
\|u^\omega-u\|_{L^q([0, T]; L^r(\R^3))}\to0\qquad\textrm{as}\;|\omega|\to\infty.
\end{equation*}
\end{theorem}
\begin{remark}
For the statement of this Theorem we restrict ourselves to the case when the power-type nonlinearity is \emph{mass-subcritical} (see \cite{T}, \cite{Caz}). Anyway the physically interesting exponent for this model, i.e. $\sigma=\frac{2}{3}$, is included in the Theorem. However, for its mathematical interest, the case of a \emph{energy-subcritical} nonlinearity will be the object of a future investigation.
\end{remark}
By means of formula \eqref{eq:psi_to_u}, the main result gives us the asymptotic behavior for $\psi$, solution of \eqref{eq:xfel}.
\begin{corollary}
Let $0<\sigma<\frac{4}{3}$, $\psi_0\in L^2(\R^3)$, $A=A^\omega(t)$ be such that
\begin{equation*}
2\int_0^tA^\omega(s)ds=b^\omega(t)=\vec e(t)f(\omega t),
\end{equation*}
as in \eqref{eq:125} and let $\psi^\omega\in\cfun(\R; L^2(\R^3))$ be the unique global solution of the Cauchy problem 
\eqref{eq:xfel}. Then for each finite time $0<T<\infty$ and for each admissible Strichartz index pair $(q, r)$ (see Section 
\ref{sect:intro}), we have
\begin{equation*}
\pare{\int_0^T\pare{\int_{\R^3}\modu{\psi^\omega(t, x)-e^{-i\int_0^t|A^\omega(s)|^2ds}u(t, x-b^\omega(t))}^rdx}^{\frac{q}{r}}dt}^
{\frac{1}{r}}=o(1),
\end{equation*}
as $|\omega|\to\infty$, where $u\in\cfun(\R; L^2(\R^3))$ is the solution to \eqref{eq:nls_aver_intro}.
\end{corollary}
In Section \ref{sect:intro} we review some results about periodic functions and weak convergence.
We also recall the Strichartz inequalities associated to the Schr\"odinger group, in the spirit of \cite{KT} (see also \cite{GV}). Such estimates will be then used in Section \ref{sect:LGWP} to perform a fixed point argument and to show the local well-posedness for the Cauchy problems \eqref{eq:nls_omega_intro} and \eqref{eq:nls_aver_intro}. By using the conservation of mass in the case of a $L^2-$subcritical power-type nonlinearity we also prove the global well-posedness (see the seminal paper by Tsutsumi \cite{Tsu}, and also the monographs \cite{Caz}, \cite{T}, \cite{LP}), by obtaining some uniform bounds for 
$\{u^\omega\}$ in $\omega$.
\newline
In Section \ref{sect:conv} we prove the main result of this paper, Theorem \ref{thm:main}. The idea for the proof is as in \cite{CS} and can be easily explained in the following way: if we consider the Duhamel's formula for equation \eqref{eq:nls_omega_intro}, 
then the oscillating potential \eqref{eq:V_omega} appears inside the time integral, thus the weak convergence for 
\eqref{eq:V_omega} can be improved to the strong one for $\{u^\omega\}$. This is indeed possible thanks to the uniform bounds in $\omega$ we have for $\{u^\omega\}$.
\section{Preliminary results and notations}\label{sect:intro}
In this Section we first recall some basic facts about weak convergence and periodic functions, which will then be extended to adapt them to our analysis. Finally, we will also give a very quick overview on dispersive estimates for the Schr\"odinger equation and on local and global analysis of its solutions.
\newline
First of all, let us recall the following theorem about weak limits of rapidly oscillating functions.
\begin{theorem}
Let $1\leq p\leq\infty$ and $f$ be a $2\pi-$periodic function in $L^p(0, 2\pi)$. Let us define
\begin{equation*}
f_n(t):=f(nt), \qquad n\in\N.
\end{equation*}
Then for $1\leq p<\infty$,
\begin{equation*}
f_n\rightharpoonup\frac{1}{2\pi}\int_0^{2\pi}f(t)dt\qquad\textrm{in}\;L^p(I),\;\textrm{for any bounded}\;\Omega\in\R,
\end{equation*}
and for $p=\infty$ we have
\begin{equation*}
f_n\rightharpoonup\frac{1}{2\pi}\int_0^{2\pi}f(t)dt\qquad\textrm{in}\;L^\infty(\R),
\end{equation*}
where the convergence is weak$-*$ in $L^\infty(\R)$.
\end{theorem}
Another very basic fact is that weak convergence is basically the convergence in average for the sequence: indeed the following theorem holds.
\begin{theorem}\label{thm:weak_conv}
Let $\{f_n\}\subset L^p(\R^N)$ be a uniformly bounded sequence in $L^p(\R^N)$. Then the following are equivalent:
\begin{enumerate}
\item 
\begin{equation*}
f_n\rightharpoonup f\qquad\textrm{in}\;L^p(\R^N);
\end{equation*}
\item \begin{equation*}
f_n\rightharpoonup f\qquad\textrm{in}\;\dist'(\R^N);
\end{equation*}
\item for each Borel set $E\subset\R^N$, $0<|E|<\infty$ we have
\begin{equation*}
\lim_{n\to\infty}\frac{1}{|E|}\int_Ef_n(x)dx=\frac{1}{|E|}\int_Ef(x)dx;
\end{equation*}
\item for each cube $E\subset\R^N$, $0<|E|<\infty$ we have
\begin{equation*}
\lim_{n\to\infty}\frac{1}{|E|}\int_Ef_n(x)dx=\frac{1}{|E|}\int_Ef(x)dx;
\end{equation*}
\end{enumerate}
\end{theorem}
\begin{remark}
The above Theorem holds for each exponent $1\leq p\leq\infty$: obviously in the $p=\infty$ case one has to change the weak convergence with the weak$-\ast$ convergence in $L^\infty$.
\end{remark}
Furthermore, the same result is also valid in the more general case of a function $g\in L^p(\R; X)$, where $X$ is an arbitrary Banach space. Clearly we are interested in the case when $X$ is a Lebesgue space $L^s(\R^d)$. Let $g\in L^p(\R; L^s(\R^d))$, such that
$g(t+2\pi, \cdot)=g(t, \cdot)$ in $L^s(\R^d)$, for each $t\in\R$, then let us define the sequence $\{g_n\}\subset L^p_tL^s_x$ in the following way:
\begin{equation*}
g_n(t, x):=g(nt, x).
\end{equation*}
Then we can prove
\begin{equation*}
g_n\rightharpoonup\frac{1}{2\pi}\int_0^{2\pi}g(t,\cdot)dt\qquad \textrm{in}\;L^p_tL^s_x.
\end{equation*}
Indeed, to prove the validity of this weak limit it suffices to prove the convergence in average on sets 
$(a, b)\times E\subset\R\times\R^d$, where $E\subset\R^d$ is a bounded Borel set in $\R^d$. Since $g\in L^p_tL^s_x$, then the function
\begin{equation*}
t\mapsto\frac{1}{|E|}\int_Eg(t, x)dx
\end{equation*}
is in $L^p(\R)$, is $2\pi-$periodic, thus it weakly converges to its average,
\begin{equation*}
\frac{1}{|E|}\int_Eg_n(\cdot, x)dx\rightharpoonup\frac{1}{2\pi}\int_0^{2\pi}\frac{1}{|E|}\int_Eg(t, x)dxdt\qquad\textrm{in}\;L^p(\R).
\end{equation*}
Hence, by the convergence in average, we have that for all $(a, b)\subset\R$
\begin{equation*}
\frac{1}{(b-a)|E|}\int_a^b\int_Eg_n(t, x)dxdt\to\frac{1}{(b-a)|E|}\int_a^b\int_E\frac{1}{2\pi}\int_0^{2\pi}g(t', x)dt'dxdt,
\end{equation*}
and this clearly means
\begin{equation*}
g_n\rightharpoonup\frac{1}{2\pi}\int_0^{2\pi}g(t, \cdot)dt\qquad\textrm{in}\;L^p_tL^s_x.
\end{equation*}
Now, let us consider a double scale function, i.e. a function depending on a \emph{slow} and a \emph{fast} variable. To best adapt the discussion below to our analysis we consider only a special class amongst those functions, namely
\begin{equation*}
\tilde g(t, \tau)=g(e(t)f(\tau)).
\end{equation*}
In our specific case $t$ will be the slow variable and $\tau$ the fast one. We assume $f$ to be $2\pi-$periodic as before, 
$e\in\cfun^\infty(\R)$ (or just regular enough) but not periodic in general, and $g$ continuous and such that, as it is defined, it lies in $L^p(\R)$. Let us define the sequence
\begin{equation*}
g_n(t):=\tilde g(t, nt)=g(e(t)f(nt)),\qquad n\in\N,
\end{equation*}
then we can show that
\begin{equation*}
g_n\rightharpoonup\langle g\rangle(t):=\frac{1}{2\pi}\int_0^{2\pi}g(e(t)f(\tau))d\tau\qquad\textrm{in}\;L^p,
\end{equation*}
where the convergence is weak$-*$ if $p=\infty$. Indeed, let us consider an interval $(a, b)\subset\R$, we have
\begin{multline*}
\int_a^bg(e(t)f(nt))dt=\frac{1}{n}\int_{na}^{nb}g\pare{e\pare{\frac{t}{n}}f(t)}dt\\
=\frac{1}{n}\sum_{k=0}^{\quadre{\frac{n(b-a)}{2\pi}}-1}\int_{na+2k\pi}^{na+2(k+1)\pi}g\pare{e\pare{\frac{t}{n}}f(t)}dt
+\frac{1}{n}\int_{na+\quadre{\frac{n(b-a)}{2\pi}}2\pi}^{nb}g\pare{e\pare{\frac{t}{n}}f(t)}dt.
\end{multline*}
Now, because of the continuity hypothesis on $g$ and $e$, for $n$ big enough we can approximate the integrals in the sum by
\begin{multline*}
\frac{1}{n}\sum_{k=0}^{\quadre{\frac{n(b-a)}{2\pi}}-1}\int_{na+2k\pi}^{2(k+1)2\pi}g\pare{e\pare{a+\frac{2k\pi}{n}}f(t)}dt
+\frac{1}{n}\int_{na+\quadre{\frac{n(b-a)}{2\pi}}2\pi}^{nb}g\pare{e\pare{\frac{t}{n}}f(t)}dt\\
=\frac{1}{n}\sum_{k=0}^{\quadre{\frac{n(b-a)}{2\pi}}-1}\int_{0}^{2\pi}g\pare{e\pare{a+\frac{2k\pi}{n}}f(t)}dt
+\frac{1}{n}\int_{na+\quadre{\frac{n(b-a)}{2\pi}}2\pi}^{nb}g\pare{e\pare{\frac{t}{n}}f(t)}dt.
\end{multline*}
Again, by the continuity hypothesis on $g$ we have, for $n$ going to infinity we can see the expression above converges to
\begin{equation*}
\int_a^b\frac{1}{2\pi}\int_0^{2\pi}g(e(t)f(\tau))d\tau dt.
\end{equation*}
Consequently, for each open set $(a, b)$ we show
\begin{equation*}
\lim_{n\to\infty}\frac{1}{b-a}\int_a^bg_n(t)dt=\frac{1}{b-a}\int_a^b\frac{1}{2\pi}\int_0^{2\pi}g(e(t)f(\tau))d\tau dt.
\end{equation*}
This clearly implies the same convergence for each Borel set $E\subset\R$:
\begin{equation*}
\lim_{n\to\infty}\frac{1}{|E|}\int_Eg_n(t)dt=\frac{1}{|E|}\int_E\frac{1}{2\pi}\int_0^{2\pi}g(e(t)f(\tau))d\tau dt,
\end{equation*}
and this, plus the uniform bound on the sequence $\{g_n\}\subset L^p(\R)$ proves that
\begin{equation*}
g_n\rightharpoonup\frac{1}{2\pi}\int_0^{2\pi}g(e(t)f(\tau))d\tau\qquad\textrm{in}\;L^p(\R).
\end{equation*}
Now, in a similar way as before, we can extend the same result to the case of $g\in L^p(\R; X)$, where $X$ is a Banach space.
\newline
In our specific case, we consider
\begin{equation}\label{eq:pot_omega}
V^\omega(t, x)=\frac{c}{|x-\vec e(t)\sin(\omega t)|}.
\end{equation}
First of all, we notice $V^\omega\in L^\infty(\R: L^{p_1}(\R^3)+L^{p_2}(\R^3))$, where $p_1, p_2$ are two Lebesgue exponents such that $p_1<3<p_2$ and are sufficiently close to 3. Indeed, let $B_t$ be the unit ball in $\R^3$ centered at the point 
$\vec e(t)\sin(\omega t)$, and let $\chi_{B_t}$ the its characteristic function. We then write 
$V^\omega=V^\omega_1+V^\omega_2:=V^\omega\chi_{B_t}+V^\omega(1-\chi_{B_t})$, and
$V^\omega_1\in L^\infty(\R; L^{p_1}(\R^3))$, $V^\omega_2\in L^\infty(\R; L^{p_2}(\R^3))$, where $p_1=\frac{3}{1+3\eps}$, 
$p_2=\frac{3}{1-3\eps}$, for some small $\eps>0$. Let us furthermore notice that the norm 
$\|V^\omega\|_{L^\infty(\R: L^{p_1}(\R^3)+L^{p_2}(\R^3))}$ does \emph{not} depend on $\omega$.
\newline
Consequently, from what we said above, we see the sequence $\{V^\omega\}$ converges weakly to the function
\begin{equation}\label{eq:aver_pot}
\langle V\rangle(t,x):=\frac{1}{2\pi}\int_0^{2\pi}\frac{c}{|x-\vec e(t)\sin(\omega\tau)|}d\tau,
\end{equation}
in $L^\infty(\R: L^{p_1}(\R^3)+L^{p_2}(\R^3))$.
\newline
Here and throughout the paper we shall set $p_1:=\frac{3}{1+3\eps}, p_2:=\frac{3}{1-3\eps}$.
\newline
Finally, let us consider also when $g$ is a regular function, in which case we have further convergence properties.
More in particular, we consider a smooth function $\zeta\in\cfun^\infty(\R^3)$, such that it is in $L^\infty(\R^3)$, together with all its derivatives. Let 
\begin{equation}\label{eq:smooth_V}
g(t, \tau, x)=\zeta(x-\vec e(t)\sin(\tau)),
\end{equation}
where $t\in[0, T]$, $\tau\in\R$, $x\in\R^3$, $0<T<\infty$ is fixed. We state a Lemma which will be useful later on this article.
\begin{lemma}\label{lemma:unif_conv}
Let $g$ be defined as in \eqref{eq:smooth_V}. Then
\begin{equation*}
\sup_{(t, x)\in[0, T]\times\R^3}\modu{
\int_0^t\pare{\zeta(x-\vec e(t')\sin(\omega t'))-\frac{1}{2\pi}\int_{-\pi}^\pi\zeta(x-\vec e(t')\sin(\tau)d\tau} dt'}\to0,
\end{equation*}
as $|\omega|\to\infty$.
\end{lemma}
\proof
Let us define
\begin{equation*}
g_l(t, x)=\frac{1}{2\pi}\int_{-\pi}^\pi g(t, \tau, x)e^{-il\tau}d\tau,\qquad l\in\Z,
\end{equation*}
so that we can write
\begin{equation*}
g(t, \tau, x)=\sum_{l\in\Z}g_l(t, x)e^{il\tau},
\end{equation*}
and we have
\begin{equation*}
\frac{1}{2\pi}\int_{-\pi}^\pi|g(t, \tau, x)|^2d\tau=\sum_{l\in\Z}|g_l(t, x)|^2.
\end{equation*}
Furthermore, we use the Fourier transform for the slow time variable. For this purpose we extend the function $g$ from 
$[0, T]\times\R\times\R^3$ to $\R\times\R\times\R^3$, such that it is smooth in $\R\times\R\times\R^3$ and it vanishes outside the slab $(-1, T+1)\times\R\times\R^3$. Thus we have
\begin{equation*}
g(t, \tau, x)=\sum_{l\in\Z}e^{il\tau}\int_{\R}e^{i\sigma t}\hat g_l(\sigma, x)d\sigma
=\sum_{l\in\Z}\int_\R \hat g_l(\sigma, x)e^{i(\sigma t+l\tau)}d\sigma.
\end{equation*}
where 
\begin{equation*}
\hat g_l(\sigma, x)=\frac{1}{2\pi}\int_\R e^{-i\sigma t}g_l(t, x)dt.
\end{equation*}
It is straightforward to see that we can write
\begin{equation*}
\int_0^t\pare{\zeta(x-\vec e(t')\sin(\omega t'))-\frac{1}{2\pi}\int_{-\pi}^\pi\zeta(x-\vec e(t')\sin(\tau))d\tau} dt'=
\int_0^tg(t', \tau, x)-g_0(t', x)dt'.
\end{equation*}
Conseuqently, to prove the statement of the Lemma, we must prove
\begin{equation*}
\int_0^t\sum_{l\neq0}\int_\R\hat g_l(\sigma, x)e^{i(\sigma+l\omega)t'}d\sigma dt'\to0,
\end{equation*}
as $|\omega|\to\infty$, uniformly in $(t, x)\in[0, T]\times\R^3$. 
\newline
Let us swap the integration order in the above expression, we then obtain
\begin{equation}\label{eq:sum_gl}
\sum_{l\neq0}\int_\R\hat g_l(\sigma, x)\frac{e^{i(\sigma+l\omega)t}-1}{\sigma+l\omega}d\sigma.
\end{equation}
Let us consider for the moment each integral in the sum,
\begin{equation*}
\int_\R\hat g_l(\sigma, x)\frac{e^{i(\sigma+l\omega)t}-1}{\sigma+l\omega}d\sigma,
\end{equation*}
without loss of generality we can consider now the case when $\omega>0$ and $l>0$. We split the above integral in two regions, inside and outside the ball centered at the origin of radius $\frac{3}{4}l\omega$. 
\begin{multline*}
\int_\R\hat g_l(\sigma, x)\frac{e^{i(\sigma+l\omega)t}-1}{\sigma+l\omega}d\sigma\\
=
\int_{\{|\sigma|\leq\frac{3}{4}l\omega\}}\hat g_l(\sigma, x)\frac{e^{i(\sigma+l\omega)t}-1}{\sigma+l\omega}d\sigma+
\int_{\{|\sigma|\geq\frac{3}{4}l\omega\}}\hat g_l(\sigma, x)\frac{e^{i(\sigma+l\omega)t}-1}{\sigma+l\omega}d\sigma.
\end{multline*}
For the first one we have that, in this region, $|\sigma+l\omega|\geq\frac{1}{4}l\omega$, hence
\begin{equation*}
\modu{\int_{\{|\sigma|\leq\frac{3}{4}l\omega\}}\hat g_l(\sigma, x)\frac{e^{i(\sigma+l\omega)t}-1}{\sigma+l\omega}d\sigma}
\leq
\frac{8}{l\omega}\int_\R|\hat g_l(\sigma, x)|d\sigma.
\end{equation*}
Now we have
\begin{multline}\label{eq:hs}
\int_\R|\hat g_l(\sigma, x)|d\sigma\leq\pare{\int_\R\frac{1}{(1+|\sigma|^2)^{\frac{1}{2}+\eps}}d\sigma}^{1/2}
\pare{\int_\R(1+|\sigma|^2)^{\frac{1}{2}+\eps}|\hat g_l(\sigma, x)|^2d\sigma}\\
\lesssim\pare{\int_\R|\langle D_t\rangle^{\frac{1+\eps}{2}}g_l(t, x)|^2dt}^{1/2}.
\end{multline}
Consequently,
\begin{equation*}
\modu{\int_{\{|\sigma|\leq\frac{3}{4}l\omega\}}\hat g_l(\sigma, x)\frac{e^{i(\sigma+l\omega)t}-1}{\sigma+l\omega}d\sigma}
\lesssim
\frac{1}{l\omega}\pare{\int_\R|\langle D_t\rangle^{\frac{1+\eps}{2}}g_l(t, x)|^2dt}^{1/2},
\end{equation*}
which tells us that the sum of those terms is $O(\frac{1}{\omega})$, as we will see more precisely later on.
On the other hand, for the second integral we use the fact that the Fourier transform of a $\cfun^\infty$ function decays faster than any polynomial, consequently we have that for each $N\in\N$, hence those integrals give us a contribution which is smaller than any power of $\omega$.
\begin{multline*}
\modu{\int_{\{|\sigma|\geq\frac{3}{4}l\omega\}}\hat g_l(\sigma, x)\frac{e^{i(\sigma+l\omega)t}-1}{\sigma+l\omega}d\sigma}
\lesssim\int_{\{|\sigma|\geq\frac{3}{4}l\omega\}}|\hat g_l(\sigma, x)|d\sigma\\
\lesssim
\frac{1}{(1+(l\omega)^2)^{N/2}}\int_\R(1+|\sigma|^2)^{N/2}|\hat g_l(\sigma, x)|d\sigma
\lesssim
\frac{1}{(1+(l\omega)^2)^{N/2}}\pare{\int_\R|\langle D_t\rangle^{\frac{N+1+\eps}{2}}g_l(t, x)|^2dt}^{1/2},
\end{multline*}
where the last inequality follows from \eqref{eq:hs}.
Thus, by taking the modulus of the sum in \eqref{eq:sum_gl}, we obtain
\begin{multline*}
\modu{\sum_{l\neq0}\int_\R\hat g_l(\sigma, x)\frac{e^{i(\sigma+l\omega)t}-1}{\sigma+l\omega}d\sigma}\\
\lesssim\frac{1}{\omega}
\sum_{l\neq0}\pare{\frac{1}{l}\pare{\int_\R|\langle D_t\rangle^{\frac{1+\eps}{2}}g_l(t, x)|^2dt}^{1/2}+
\frac{1}{(1+(l\omega)^2)^{1/2}}\pare{\int_\R|\langle D_t\rangle^{\frac{3+\eps}{2}}g_l(t, x)|^2dt}^{1/2}}.
\end{multline*}
By Cauchy-Schwartz inequality in the sum above we get
\begin{multline*}
\modu{\sum_{l\neq0}\int_\R\hat g_l(x, \sigma)\frac{e^{i(\sigma+l\omega)t}-1}{\sigma+l\omega}d\sigma}
\lesssim
\frac{1}{\omega}\pare{\sum_{l\neq0}\frac{1}{l^2}}^{1/2}
\pare{\sum_{l\neq0}\int_\R|\langle D_t\rangle^{\frac{3+\eps}{2}}g(x, t)|^2dt}^{1/2}\\
\lesssim
\frac{1}{\omega}\pare{\int_{-\pi}^\pi\int_\R|\langle D_t\rangle^{\frac{3+\eps}{2}}g(x, t, \tau)|^2dtd\tau}^{1/2}.
\end{multline*}
Hence we can conclude that
\begin{multline*}
\sup_{(t, x)\in[0, T]\times\R^3}
\modu{\int_0^t\pare{\zeta(x-\vec e(t')\sin(\omega t'))-\frac{1}{2\pi}\int_{-\pi}^\pi \zeta(x-\vec e(t')\sin(\tau))d\tau} dt'}\\
\lesssim
\frac{1}{\omega}\sup_{x\in\R^3}
\pare{\int_{-\pi}^\pi\int_\R|\langle D_t\rangle^{\frac{3+\eps}{2}}g(t, \tau, x)|^2dtd\tau}^{1/2}\\
=
\frac{1}{\omega}
\pare{\sup_{x\in\R^3}\int_{-\pi}^\pi\int_{-1}^{T+1}|\langle D_t\rangle^{\frac{3+\eps}{2}}g(t, \tau, x)|^2dtd\tau}^{1/2}\\
\lesssim\frac{2\pi (T+2)}{\omega}\|\langle D\rangle^{\frac{3+\eps}{2}}\zeta\|_{L^\infty(\R^3)},
\end{multline*}
which proves the Lemma.
\qed
\subsection{Review of Strichartz estimates and Local and Global Theory for Schr\"odinger Equations}\label{sub:strich}
In this subsection we quickly review some basic facts about dispersive estimates for the Schr\"odinger equation and their application to local and global existence theory of solutions to nonlinear Schr\"odinger equations.
\newline
Let $U(t):=e^{it\Delta}$ denote the free Schr\"odinger group, i.e. if $u$ is solution to
\begin{equation*}
\left\{\begin{array}{l}
i\d_tu=-\Delta u\\
u(0)=u_0,
\end{array}\right.
\end{equation*}
then $u(t)=U(t)u_0$.
\begin{definition}\label{defn:admiss}
We say $(q, r)$ is an \emph{admissible pair} of exponents if $2\leq q\leq\infty, 2\leq r\leq6$, and 
\begin{equation*}
\frac{1}{q}=\frac{3}{2}\pare{\frac{1}{2}-\frac{1}{r}}.
\end{equation*}
\end{definition}
Here and throughout the paper we will say $\|\cdot\|_{L^q_tL^r_x}$ is a \emph{Strichartz norm} to mean that it is a norm taken in a space such that $(q, r)$ is an admissible pair in the sense of Definition \ref{defn:admiss}. We will also use the notation
\begin{equation}\label{eq:strich}
\|f\|_{S(I)}:=\sup_{(q, r)}\|f\|_{L^q_tL^r(I\times\R^3)},
\end{equation}
where the $\sup$ is taken over all admissible pairs $(q, r)$.
\newline
Furthermore, we also need
\begin{equation*}
\|f\|_{S^1(I)}:=\|f\|_{S(I)}+\|\nabla f\|_{S(I)}.
\end{equation*}
\newline
Now we write the Strichartz estimates we will need in our paper. Such estimates  go back to Strichartz \cite{Str} which proved them in a particular case for the wave equation, then Ginibre, Velo \cite{GV}, through a $TT^\ast$ argument, extended the result and finally Keel, Tao \cite{KT} proved the \emph{endpoint estimate} and general dispersive estimates in an abstract setup. Such estimates hold for a general dispersive equation in arbitrary space dimensions, but for our study we will use (and state) them only for the Schr\"odinger equation in $\R^3$.
\begin{theorem}[Keel, Tao \cite{KT}]
Let $(q, r), (\tilde q, \tilde r)$ be two arbitrary admissible pairs and let $U(\cdot)$ be the free Schr\"odinger group. Then
\begin{align*}
\|U(t)f\|_{L^q_tL^r_x}&\lesssim\|f\|_{L^2}\\
\|\int_0^tU(t-s)F(s)ds\|_{L^q_tL^r_x}&\lesssim\|F\|_{L^{\tilde q'}_tL^{\tilde r'}_x}\\
\|\int U(-s)F(s)ds\|_{L^2}&\lesssim\|F\|_{L^{\tilde q'}_tL^{\tilde r'}_x}.\\
\end{align*}
\end{theorem}
Strichartz estimates are very useful to prove existence of local solutions to nonlinear Schr\"odinger equations through a fixed point argument.
\newline
Indeed, let us consider the following Schr\"odinger equation
\begin{equation*}
\left\{\begin{array}{l}
i\d_tu=-\Delta u+F_1+\dotsc+F_M\\
u(0)=u_0,
\end{array}\right.
\end{equation*}
for some functions $F_1, \dotsc, F_M$, then by Duhamel's formula we can write
\begin{equation*}
u(t)=U(t)u_0-i\int_0^tU(t-s)(F_1+\dotsc F_M)(s)ds.
\end{equation*}
Then, by applying Strichartz estimates to the above formula we obtain
\begin{equation*}
\|u\|_{L^q_tL^r_x}\lesssim\|u_0\|_{L^2}+\|F_1\|_{L^{q_1'}_tL^{r_1'}_x}+\dotsc+\|F_M\|_{L^{q_M'}_tL^{r_M'}_x},
\end{equation*}
where $(q, r), (q_1, r_1), \dotsc, (q_M, r_M)$ are admissible pairs. Hence it is clear that we can bound each term $F_j$ in an arbitrary dual Strichartz space $L^{q_j'}_tL^{r_j'}_x$. If now the $F_j$'s are different nonlinearities, we further estimate each term 
$\|F_j\|_{L^{q_j'}_tL^{r_j'}_x}$ to close the fixed point argument. The reader should see the monographs \cite{Caz}, \cite{T} and references therein for details.
\section{Local and Global Well-Posedness}\label{sect:LGWP}
In this Section we state the local and global well-posedness results we have for equations \eqref{eq:nls_omega_intro} and 
\eqref{eq:nls_aver_intro}. We first prove a local well-posedness result for \eqref{eq:nls_omega_intro} in the space of energy (i.e. $H^1(\R^3)$), then global well-posedness in the space of mass (i.e. $L^2(\R^3)$). Furthermore we show some uniform bounds on the $S(0, T)$ (see \eqref{eq:strich}) norm of solutions to equation \eqref{eq:nls_omega_intro}. 
As we already mentioned in Section \ref{sect:intro}, both $V^\omega$ and $\langle V\rangle$ belong to the space 
$L^\infty(\R: L^{p_1}(\R^3)+L^{p_2}(\R^3))$, where $p_1=\frac{3}{1+3\eps}, p_2=\frac{3}{1-3\eps}$, with $\eps>0$ sufficiently small, and the norm $\|V^\omega\|_{L^\infty(\R: L^{p_1}(\R^3)+L^{p_2}(\R^3))}$ does \emph{not} depend on $\omega$.
 However, here we want to show a local well-posedness theory in $H^1$, hence we also need to estimate their gradient. We see that
\begin{equation*}
\nabla V^\omega, \nabla\langle V\rangle\in L^\infty(\R; L^{\frac{3}{2+3\eps}}(\R^3)+L^{\frac{3}{2-3\eps}}(\R^3)).
\end{equation*}
Indeed, let us consider again the characteristic function $\chi_{B_t}$ of the unit ball $B_t$ centered at 
$\vec e(t)\sin(\omega t)\in\R^3$, then 
$\chi_{B_t}\nabla V^\omega\in L^\infty(\R; L^{\frac{3}{2+3\eps}}(\R^3))$ and 
$(1-\chi_{B_t})\nabla V^\omega\in L^\infty(\R; L^{\frac{3}{2-3\eps}}(\R^3))$. Again, the norm 
$\|\nabla V^\omega\|_{L^\infty(\R; L^{\frac{3}{2+3\eps}}(\R^3)+L^{\frac{3}{2-3\eps}}(\R^3))}$ does \emph{not} depend on 
$\omega$. Consequently, by the weak convergence, we can say that $\langle V\rangle$ is in 
$L^\infty(\R; L^{\frac{3}{2+3\eps}}(\R^3)+L^{\frac{3}{2-3\eps}}(\R^3))$, too.
\newline
Consequently here we investigate the Cauchy problem
\begin{equation}\label{eq:nls}
\left\{\begin{array}{l}
i\d_tu=-\Delta u+Vu+C_1(|\cdot|^{-1}\ast|u|^2)u-a|u|^\sigma u\\
u(0)=u_0,
\end{array}\right.
\end{equation}
where $V$ is a general potential
 such that $V=V_1+V_2$, where 
$V_1\in L^\infty(\R;L^{p_1}(\R^3))$ and $V_2\in L^\infty(\R;L^{p_2}(\R^3))$, and 
$\nabla V=\nabla V_1+\nabla V_2$, where
$\nabla V_1\in L^\infty(\R;L^{\frac{3}{2+3\eps}}(\R^3))$ and
$\nabla V_2\in L^\infty(\R;L^{\frac{3}{2-3\eps}}(\R^3))$.
In this way we have the well-posedness results below apply both to equation \eqref{eq:nls_omega_intro} and \eqref{eq:nls_aver_intro}.
\begin{theorem}\label{thm:LWP}
Assume $0<\sigma<4$, $u_0\in H^1(\R^3)$. Then, there exists a unique local solution to \eqref{eq:nls}. Furthermore, we have $u\in S^1(0, T)$.
\end{theorem}
The proof of the Theorem above is standard, based on a fixed point argument, see for example \cite{Caz}, \cite{T}, \cite{GV2}, \cite{K}. Nevertheless, for the sake of clarity and completeness, we sketch the main steps.
\newline
For this purpose. let us first write two techincal lemmas which will be used in what follows.
\begin{lemma}\label{lemma:hartree_power}
Let $k=0, 1$. Then
\begin{equation}\label{eq:i}
\|\nabla^k(|\cdot|^{-1}\ast|f|^2)f\|_{L^1_tL^2_x([0, T]\times\R^3)}
\lesssim T^{\frac{1}{2}}\|f\|_{L^6_tL^{18/7}_x}\|\nabla^kf\|_{L^6_tL^{18/7}_x};
\end{equation}
\begin{equation}\label{eq:ii}
\|\nabla^k(|f|^\sigma f)\|_{L^{\frac{4(\sigma+2)}{\sigma+8}}_tL^{\frac{\sigma+2}{\sigma+1}}_x}
\lesssim T^{\frac{2(\sigma+2)}{4-\sigma}}\|f\|_{L^\infty_tH^1_x}^\sigma
\|\nabla^kf\|_{L^{\frac{4(\sigma+2)}{\sigma+8}}_tL^{\sigma+2}_x};
\end{equation}
\begin{multline}\label{eq:iii}
\|(|\cdot|^{-1}\ast|f|^2)f-(|\cdot|^{-1}\ast|g|^2)g\|_{L^1_tL^2_x([0, T]\times\R^3)}\\
\lesssim T^{\frac{1}{2}}\pare{\|f\|_{L^6_tL^{18/7}_x}^2+\|g\|_{L^6_tL^{18/7}_x}^2}\|f-g\|_{L^6_tL^{18/7}_x};
\end{multline}
\begin{multline}\label{eq:iv}
\||f|^\sigma f-|g|^\sigma g\|_{L^{\frac{4(\sigma+2)}{\sigma+8}}_tL^{\frac{\sigma+2}{\sigma+1}}_x}\\
\lesssim T^{\frac{2(\sigma+2)}{4-\sigma}}
\pare{\|f\|_{L^\infty_tH^1_x}^\sigma+\|g\|_{L^\infty_tH^1_x}^\sigma}
\|f-g\|_{L^{\frac{4(\sigma+2)}{\sigma+8}}_tL^{\sigma+2}_x};
\end{multline}
\end{lemma}
\proof
\eqref{eq:i} is a consequence of the following inequality
\begin{equation}\label{eq:ast}
\|(|\cdot|^{-1}\ast(f_1f_2))f_3\|_{L^1_tL^2_x}\lesssim
T^{\frac{1}{2}}\prod_{i=1}^3\|f_i\|_{L^6_tL^{18/7}_x},
\end{equation}
which can be proved by using H\"older's inequality and Hardy-Littlewood-Sobolev inequality.
\newline
In a similar way by H\"older's inequality we get
\begin{equation}\label{eq:astast}
\||f_1|^\sigma f_2\|_{L^{\frac{4(\sigma+2)}{\sigma+8}}_tL^{\frac{\sigma+2}{\sigma+1}}_x}
\lesssim T^{\frac{2(\sigma+2)}{4-\sigma}}\|f_1\|_{L^\infty_tL^{\sigma+2}_x}^\sigma
\|f_2\|_{L^{\frac{4(\sigma+2)}{\sigma+8}}_tL^{\sigma+2}_x},
\end{equation}
and then, since $0<\sigma<4$, we can use Sobolev embedding to show that 
\begin{equation*}
\||f_1|^\sigma f_2\|_{L^{\frac{4(\sigma+2)}{\sigma+8}}_tL^{\frac{\sigma+2}{\sigma+1}}_x}
\lesssim T^{\frac{2(\sigma+2)}{4-\sigma}}\|f_1\|_{L^\infty_tH^1_x}^\sigma
\|f_2\|_{L^{\frac{4(\sigma+2)}{\sigma+8}}_tL^{\sigma+2}_x}.
\end{equation*}
Now, \eqref{eq:ii} follows from the fact that $|\nabla(|f|^\sigma f)|\lesssim|f|^\sigma|\nabla f|$.
\newline
The proof of inequalitites \eqref{eq:iii} and \eqref{eq:iv} are similar to \eqref{eq:i} and \eqref{eq:ii}, respectively. In particular, by using some algebra and \eqref{eq:ast} yields \eqref{eq:iii}, whereas \eqref{eq:iv} follows from \eqref{eq:astast} and
\begin{equation}
\modu{|f|^\sigma f-|g|^\sigma g}\lesssim\pare{|f|^\sigma+|g|^\sigma}|f-g|.
\end{equation}
\qed\newline
The second technical lemma estimates the terms with the Coulomb potentials, both in \eqref{eq:nls_omega_intro} and in 
\eqref{eq:nls_aver_intro}. 
\begin{lemma}
\begin{equation}\label{eq:uno}
\|\nabla Vf\|_{L^2_tL^{6/5}_x([0, T]\times\R^3}
\leq T^{\frac{3\eps}{2}}\|\nabla V_1\|_{L^\infty_tL^{\frac{3}{2+3\eps}}_x}\|f\|_{L^{\frac{2}{1-3\eps}}_tL^{\frac{6}{1+6\eps}}};
\end{equation}
\begin{equation}\label{eq:due}
\|V_1(1+|\nabla|)f\|_{L^2_tL^{6/5}_x}
\leq T^{\frac{1-3\eps}{2}}\|V_1\|_{L^\infty_tL^{\frac{3}{1+3\eps}}_x}\|(1+|\nabla|)f\|_{L^{\frac{2}{3\eps}}_tL^{\frac{2}{1-2\eps}}_x};
\end{equation}
\begin{equation}\label{eq:tre}
\|V_2(1+|\nabla|)f\|_{L^1_tL^2_x([0, T]\times\R^3}\leq T^{\frac{1+3\eps}{2}}
\|V_2\|_{L^\infty_tL^{\frac{3}{1-3\eps}}_x}\|(1+|\nabla|)f\|_{L^{\frac{2}{1-3\eps}}_tL^{\frac{6}{1+6\eps}}_x}.
\end{equation}
\end{lemma}
\begin{remark}
Let us notice that the pair of exponents 
$(6, \frac{18}{7})$, $\pare{\frac{4(\sigma+2)}{3\sigma}, \sigma+2}$, $\pare{\frac{2}{1-3\eps}, \frac{6}{1+6\eps}}$,
$\pare{\frac{2}{3\eps}, \frac{2}{1-2\eps}}$, are all Schr\"odinger admissible, thus the norms in those spaces are all bounded by the $S(I)$ norm.
\end{remark}
We can now sketch the proof of Theorem \ref{thm:LWP}
\proof
Let $u_0\in H^1(\R^3)$ be given. By the Duhamel's formula we have
\begin{equation*}
u(t)=e^{it\Delta}u_0-i\int_0^te^{i(t-s)\Delta}\quadre{Vu+C_1(|\cdot|^{-1}\ast|u|^2)u-a|u|^\sigma u}(s)ds.
\end{equation*}
We want to show that, for $0<T$ sufficiently small,
\begin{equation*}
\Phi(w):=e^{it\Delta}u_0-i\int_0^te^{i(t-s)\Delta}\quadre{Vu+C_1(|\cdot|^{-1}\ast|u|^2)u-a|u|^\sigma u}(s)ds
\end{equation*}
maps a ball $B_R\subset S^1(0, T)$ (the radius $R$ will be chosen later) into itself, and that in this ball $\Phi$ is a contraction in the $S(0, T)$ metric.
\newline
By Strichartz estimates we have
\begin{multline*}
\|\Phi(w)\|_{S^1(0, T)}\lesssim\|u_0\|_{H^1}
+\|\nabla Vw\|_{L^2_tL^{6/5}_x}
+\|V_1(1+|\nabla|)w\|_{L^2_tL^{6/5}_x}
+\|V_2(1+|\nabla|)w\|_{L^1_tL^2_x}\\
+\|(1+|\nabla|)\pare{(|\cdot|^{-1}\ast|w|^2)w}\|_{L^1_tL^2_x}
+\||w|^\sigma(1+|\nabla|)w\|_{L^{\frac{4(\sigma+2)}{\sigma+8}}_tL^{\frac{\sigma+2}{\sigma+1}}_x}.
\end{multline*}
Now we use inequalities \eqref{eq:uno}, \eqref{eq:due}, \eqref{eq:tre}, \eqref{eq:i}, \eqref{eq:ii} to obtain
\begin{align*}
\|\Phi(w)\|_{S^1(0, T)}\lesssim\|u_0\|_{H^1}
&+T^{\frac{3\eps}{2}}\|\nabla V_1\|_{L^\infty_tL^{\frac{3}{2+3\eps}}_x}\|w\|_{S^1}\\
&+T^{\frac{1-3\eps}{2}}\|\nabla V_2\|_{L^\infty_tL^{\frac{3}{2-3\eps}}_x}\|w\|_{S^1}\\
&+ T^{\frac{1-3\eps}{2}}\|V_1\|_{L^\infty_tL^{\frac{3}{1+3\eps}}_x}\|w\|_{S^1}\\
&+T^{\frac{1+3\eps}{2}}\|V_2\|_{L^\infty_tL^{\frac{3}{1-3\eps}}_x}\|w\|_{S^1}\\
&+T^{\frac{1}{2}}\|w\|_{S^1}^3\\
&+T^{\frac{2(\sigma+2)}{4-\sigma}}\|w\|_{S^1}^{\sigma+1}.
\end{align*}
Thus, if we take $0<T\leq1$ sufficiently small, we have $\Phi:B_r\to B_R\subset S^1(0, T)$, for some radius depending on 
$\|u_o\|_{H^1}$. Furthermore, by using \eqref{eq:iii} and \eqref{eq:iv}, we also have
\begin{align*}
\|\Phi(w_1)-\Phi(w_2)\|_{S^1(0, T)}\lesssim
&+T^{\frac{3\eps}{2}}\|\nabla V_1\|_{L^\infty_tL^{\frac{3}{2+3\eps}}_x}\|w_1-w_2\|_{S(0, T)}\\
&+T^{\frac{1-3\eps}{2}}\|\nabla V_2\|_{L^\infty_tL^{\frac{3}{2-3\eps}}_x}\|w_1-w_2\|_{S(0, T)}\\
&+ T^{\frac{1-3\eps}{2}}\|V_1\|_{L^\infty_tL^{\frac{3}{1+3\eps}}_x}\|w_1-w_2\|_{S(0, T)}\\
&+T^{\frac{1+3\eps}{2}}\|V_2\|_{L^\infty_tL^{\frac{3}{1-3\eps}}_x}\|w_1-w_2\|_{S(0, T)}\\
&+T^{\frac{1}{2}}\pare{\|w_1\|_{S^1}^2+|w_2\|_{S^1}^2}\|w_1-w_2\|_{S(0, T)}\\
&+T^{\frac{2(\sigma+2)}{4-\sigma}}\pare{\|w_1\|_{S^1}^\sigma+\|w_2\|_{S^1}^\sigma}\|w_1-w_2\|_{S(0, T)}.
\end{align*}
Once again, if we take $0<T\leq1$ small enough, then we have $\Phi:B_R\to B_R$ is a contraction in the $S(0, T)$ metric. Thus there exists a fixed point for $\Phi$ which is hence a local solution for \eqref{eq:nls}.
\newline\qed\newline
Next Theorem deals with the well-posedness issue in $L^2(\R^3)$ for \eqref{eq:nls_omega_intro} and \eqref{eq:nls_aver_intro}. We show that when the power $\sigma$ is \emph{mass-subcritical}, i.e. $0<\sigma<\frac{4}{3}$, then for any initial datum in $L^2$ we have a \emph{global} solution. We also stress here that the uniform bound we obtain for the $S(0, T)$ norm of the solution does \emph{not} depend on $\omega$.
\begin{theorem}\label{thm:mass_subcrit}
Assume $0<\sigma<4/3$ and consider $u_0\in L^2(\R^3)$. The solution for \eqref{eq:nls} is global, 
$u\in C(\R: L^2(\R^3))$. Furthermore for each finite time $0<T<\infty$ we have
\begin{equation}\label{eq:unif_bounds}
\|u\|_{S(0, T)}\leq C(\|u_0\|_{L^2(\R^3)}, T),
\end{equation}
where the constant in the right hand side depends only on the $L^2$-norm of the initial datum and the time $T$, in particular it does not depend on $\omega$.
\end{theorem}
\proof
The proof works as for Theorem \ref{thm:LWP} at a local level and then we use the conservation of mass to extend the solution globally.
\newline
Let us consider the Duhamel's formula
\begin{equation*}
u(t)=e^{it\Delta}u_0-i\int_0^te^{i(t-s)\Delta}\pare{Vu+C_1(|\cdot|^{-1}\ast|u|^2)u-a|u|^\sigma u}(s)ds,
\end{equation*}
then by applying the Strichartz estimates to the expression above in the time interval $[0, T_1]$, we get
\begin{multline*}
\|u\|_{S(0, T_1)}\lesssim\|u_0\|_{L^2(\R^3)}
+T^{\frac{1-3\eps}{2}}\|V_1\|_{L^\infty_tL^{\frac{3}{1+3\eps}}_x}\|u\|_{S(0, T_1)}
+T^{\frac{1+3\eps}{2}}\|V_2\|_{L^\infty_tL^{\frac{3}{1-3\eps}}_x}\|u\|_{S(0, T_1)}\\
+T_1^{1/2}\|u\|_{S(0, T_1)}^3+T_1^{\frac{4-3\sigma}{4}}\|u\|_{S(0, T_1)}^{\sigma+1}.
\end{multline*}
Now we can see that if we choose $T_1=T_1(\|u_0\|_{L^2})$ small enough, then by a standard fixed point argument we have
\begin{equation}\label{eq:first_step}
\|u\|_{S(0, T_1)}\leq C\|u_0\|_{L^2}.
\end{equation}
Furthermore, the total mass is conserved at all times, $\|u(t)\|_{L^2}=\|u_0\|_{L^2}$. 
Thus we can repeat the argument to continue the solution also in the time interval $[T_1, 2T_1]$, and so on. Consequently the solution is global, and for any finite time $0<T<\infty$, we consider 
$[0, T]\subset[0, T_1]\cup\dotsc\cup[(N-1)T_1, NT_1]$, $N=\quadre{\frac{T}{T_1}}+1$, where here $[\cdot]$ denotes the integer part of the number. Thus
\begin{equation*}
\|u\|_{S(0, T)}\leq CN\|u_0\|_{L^2(\R^3)},
\end{equation*}
where $C$ is the constant appearing in \eqref{eq:first_step}. Consequently
\begin{equation*}
\|u\|_{S(0, T)}\leq C(\|u_0\|_{L^2(\R^3)}, T),
\end{equation*}
for each finite time $0<T<\infty$.
\qed
\begin{remark}
Here is a couple of remarks about the Theorem above.
\begin{itemize}
\item
Regarding the case when the power nonlinearity is mass-supercritical, energy-subcritical (i.e. $\frac{4}{3}<\sigma<4$), we cannot establish a global well-posedness result in $H^1(\R^3)$, not even in the defocusing case (i.e. $a<0$), because the energy is not conserved in our model, and the bound on the $H^1$-norm of the solution $u(t)$ at time $t$ in general would depend on $\omega$.
\item
On the other hand, the time dependence of the constant in \eqref{eq:unif_bounds} is unavoidable, because of the mass-subcriticality of the power-type nonlinearity. This is indeed what also happens for the usual NLS equation (see \cite{T} for instance).
\end{itemize}
\end{remark}
\section{Convergence Result}\label{sect:conv}
In this Section we will prove the main result of this paper, namely Theorem \ref{thm:main}. As already introduced in Section \ref{sect:first}, we want to show the convergence of solutions for
\begin{equation}\label{eq:nls_omega}
\left\{\begin{array}{l}
i\d_tu^\omega=-\Delta u^\omega+V^\omega u^\omega+C_1(|\cdot|^{-1}\ast|u^\omega|^2)u^\omega-a|u^\omega|^\sigma u^\omega\\
u^\omega(0)=u_0,
\end{array}\right.
\end{equation}
where $V^\omega$ is defined in \eqref{eq:pot_omega}, to solutions of
\begin{equation}\label{eq:nls_aver}
\left\{\begin{array}{l}
i\d_tu=-\Delta u+\langle V\rangle u+C_1(|\cdot|^{-1}\ast|u|^2)u-a|u|^\sigma u\\
u(0)=u_0,
\end{array}\right.
\end{equation}
where the averaged potential $\langle V\rangle$ has been defined in \eqref{eq:aver_pot}.
\newline
Let us recall the definition of the Strichartz space, already given in Section \ref{sect:intro}, which is
\begin{equation*}
\|f\|_{S(0, T)}:=\sup_{(q, r)}\|f\|_{L^q_tL^r_x([0, T]\times\R^3)},
\end{equation*}
where the $\sup$ is taken over all admissible pairs $(q, r)$.
\newline
The key point of the convergence result stated in Theorem \ref{thm:main} is the Lemma below: indeed the weak convergence of $V^\omega$ towards $\langle V\rangle$ improves to strong convergence for $u^\omega$ towards $u$ because, by considering the difference of the Duhamel's formulas for \eqref{eq:nls_omega} and \eqref{eq:nls_aver}, the term $V^\omega-\langle V\rangle$ appears inside the time integral, and thus the convergence in average for the oscillating potential yields the strong convergence for the solutions. This fact, together with the uniform bounds \eqref{eq:unif_bounds}, provides us the right convergence result. 
A similar result is considered also in \cite{CS}, where the authors study the solutions of a nonlinear Schr\"odinger equation with an oscillating nonlinearity and their asymptotic behavior when the oscillations are increasing more and more. The Lemma above is heavily inspired by Proposition 2.5 in \cite{CS}.
\begin{lemma}\label{lemma:main}
Let $0<T<\infty$, $f\in S(0, T)$, and let 
$V^\omega, \langle V\rangle$ be defined in \eqref{eq:pot_omega}, \eqref{eq:aver_pot}. Then
\begin{equation}\label{eq:conv}
\|\int_0^te^{i(t-s)\Delta}\pare{(V^\omega-\langle V\rangle)f}(s)ds\|_{S(0, T)}\to 0,
\end{equation}
as $|\omega|\to\infty$.
\end{lemma}
\proof
First of all let us point out that the norm appearing in \eqref{eq:conv} is uniformly bounded. Indeed by using Strichartz estimates we have
\begin{multline}\label{eq:a_priori_bound}
\|\int_0^te^{i(t-s)\Delta}(Vf)(s)ds\|_{S(0, T)}
\lesssim\|V_1f\|_{L^2_tL^{6/5}_x}+\|V_2f\|_{L^1_tL^2_x}\\
\lesssim T^{\frac{1-3\eps}{2}}\|V_1\|_{L^\infty_tL^{p_1}}\|f\|_{L^{\frac{2}{3\eps}}_tL^{\frac{2}{1-2\eps}}_x}
+T^{\frac{1+3\eps}{2}}\|V_2\|_{L^\infty_tL^{p_2}_x}\|f\|_{L^{\frac{2}{1-3\eps}}_tL^{\frac{6}{1+6\eps}}_x}\\
\lesssim T^\alpha\|V\|_{L^\infty([0, T]:L^{p_1}(\R^3)+L^{p_2}(\R^3))}\|f\|_{S(0, T)}.
\end{multline}
where $V$ can be either $V^\omega$ or $\langle V\rangle$, and $\alpha>0$. For the integral in \eqref{eq:conv} we have the following uniform bound
\begin{multline*}
\|\int_0^te^{i(t-s)\Delta}\quadre{(V^\omega-\langle V\rangle)f}(s)ds\|_{S(0, T)}
\lesssim\\ T^\alpha(\|V^\omega\|_{L^\infty([0, T]:L^{p_1}(\R^3)+L^{p_2}(\R^3))}
+\|\langle V\rangle\|_{L^\infty([0, T]:L^{p_1}(\R^3)+L^{p_2}(\R^3))})
\|f\|_{S(0, T)}.
\end{multline*}
Thanks to this a priori bound, by using a standard density argument it suffices to prove \eqref{eq:conv} only for 
$V^\omega, \langle V\rangle, f\in\cfun^\infty_0(\R\times\R^3)$. For the sake fo clarity we explain the last statement more in detalis. Let
$\{f_n\}\subset\cfun^\infty_0(\R^{1+3})$ be a sequence of compactly supported smooth functions such that $f_n\to f$ in $S(0, T)$.
Furthermore, set
\begin{equation*}
V^\omega_n:=\phi_{\frac{1}{n}}\ast V^\omega,
\end{equation*}
where $\phi_{\frac{1}{n}}(x)$ is a Gaussian with variance equal to $\frac{1}{n}$, and the convolution is only in space.
The following properties for $\{V^\omega_n\}$ hold true:
\begin{itemize}
\item $V^\omega_n\in\cfun^\infty(\R\times\R^3)$ and $D^\alpha V^\omega_n\in L^\infty(\R^3)$, for each multlindex $\alpha$;
\item $V^\omega_n\to V^\omega$ in $L^\infty(\R;L^{p_1}(\R^3)+L^{p_2}(\R^3))$, as $n\to\infty$;
\item $V^\omega_n\rightharpoonup\langle V\rangle_n$, where $\langle V\rangle_n=\phi_{\frac{1}{n}}\ast\langle V\rangle$ is the regularisation of the averaged potential.
\end{itemize}
Clearly we also have $\langle V\rangle_n\in\cfun(\R\times\R^3)$, $\langle V\rangle_n\to\langle V\rangle$ in $L^\infty(\R;L^{p_1}(\R^3)+L^{p_2}(\R^3))$, as $n\to\infty$.
We now consider the integral in \eqref{eq:conv}, we split it into four parts
\begin{align*}
\|\int_0^te^{i(t-s)\Delta}&\quadre{(V^\omega-\langle V\rangle)f}(s)ds\|_{S(0, T)}\\
\leq&\|\int_0^te^{i(t-s)\Delta}\quadre{(V^\omega-\langle V\rangle)(f-f_n)}(s)ds\|_{S(0, T)}\\
&+\|\int_0^te^{i(t-s)\Delta}\quadre{(V^\omega-V^\omega_n)f_n}(s)ds\|_{S(0, T)}\\
&+\|\int_0^te^{i(t-s)\Delta}\quadre{(\langle V\rangle_n-\langle V\rangle)f_n}(s)ds\|_{S(0, T)}\\
&+\|\int_0^te^{i(t-s)\Delta}\quadre{(V^\omega_n-\langle V\rangle_n)f_n}(s)ds\|_{S(0, T)}\\
=:&I_1+I_2+I_3+I_4.
\end{align*}
Let us apply the estimate \eqref{eq:a_priori_bound} to the terms $I_1, I_2, I_3$, we then obtain
\begin{align*}
I_1\leq& T^\alpha(\|V^\omega\|_{L^\infty([0, T]:L^{p_1}(\R^3)+L^{p_2}(\R^3))}
+\|\langle V\rangle\|_{L^\infty([0, T]:L^{p_1}(\R^3)+L^{p_2}(\R^3))})
\|f-f_n\|_{S(0, T)},\\
I_2\leq &T^\alpha\|V^\omega-V^\omega_n\|_{L^\infty([0, T]:L^{p_1}(\R^3)+L^{p_2}(\R^3))}\|f_n\|_{S(0, T)},\\
I_3\leq &T^\alpha\|\langle V\rangle_n-\langle V\rangle\|_{L^\infty([0, T]:L^{p_1}(\R^3)+L^{p_2}(\R^3))}\|f_n\|_{S(0, T)}.
\end{align*}
Thus it only remains to estimate the integral
\begin{equation}\label{eq:conv_reg}
\int_0^te^{i(t-s)\Delta}\quadre{(V^\omega_n-\langle V\rangle_n)f_n}(s)ds.
\end{equation}
By integrating by parts we have
\begin{equation*}
\int_0^te^{i(t-s)\Delta}\curly{\frac{d}{ds}\quadre{\int_0^sV^\omega_n-\langle V_n\rangle ds'f_n(s)}
-\int_0^sV^\omega_n-\langle V_n\rangle ds'\d_sf_n(s)}ds,
\end{equation*}
and then again we have
\begin{multline*}
\int_0^t\frac{d}{ds}\quadre{e^{i(t-s)\Delta}\pare{\int_0^sV^\omega_n-\langle V_n\rangle ds'f_n(s)}}ds\\
+i\int_0^te^{i(t-s)\Delta}\Delta\pare{\int_0^sV^\omega_n-\langle V_n\rangle ds'f_n(s)}ds
-\int_0^te^{i(t-s)\Delta}\pare{\int_0^sV^\omega_n-\langle V_n\rangle ds'\d_sf_n(s)}ds.
\end{multline*}
Hence, we get the integral in \eqref{eq:conv_reg} equals
\begin{align*}
\int_0^t&V^\omega_n-\langle V_n\rangle ds\;f_n(t)\\
&+i\int_0^te^{i(t-s)\Delta}\Big[\int_0^s\Delta V^\omega_n-\Delta\langle V_n\rangle ds'\;f_n(s)
+2\int_0^s\nabla V^\omega_n-\nabla\langle V_n\rangle ds'\cdot\nabla f_n\\
&\qquad+\int_0^sV^\omega_n-\langle V_n\rangle ds'\Delta f_n(s)\Big]ds\\
&-\int_0^te^{i(t-s)\Delta}\pare{\int_0^sV^\omega_n-\langle V_n\rangle ds'\d_sf_n(s)}ds.
\end{align*}
Now we estimate each term in the expression above in the space $S(0, T)$. By Strichartz estimates the first term is bounded by
\begin{equation*}
\|\int_0^tV^\omega_n-\langle V\rangle ds f_n(t)\|_{S(0, T)}
\leq\|\int_0^t\pare{V^\omega_n-\langle V_n\rangle}ds\|_{L^\infty_{t, x}([0, T]\times\R^3)}\|f_n\|_{S(0, T)}.
\end{equation*}
The other terms are estimated similarly, let us consider the next one for example
\begin{multline*}
\|\int_0^te^{i(t-s)\Delta}\Big[\int_0^s\Delta V^\omega-\langle \Delta V_n\rangle ds'f_n(s)\Big]ds\|_{S(0, T)}\\
\lesssim\|\int_0^t\pare{\Delta V^\omega_n-\langle\Delta V_n\rangle}ds\|_{L^\infty_{t, x}([0, T]\times\R^3)}
\|f_n\|_{L^{\tilde q'}_tL^{\tilde r'}_x},
\end{multline*}
for some admissible pair $(\tilde q, \tilde r)$.
\newline
Consequently, by putting everything together we obtain
\begin{align*}
\|\int_0^te^{i(t-s)\Delta}&\quadre{(V^\omega_n-\langle V\rangle_n)f_n}(s)ds\|_{S(0, T)}\\
\lesssim&\|\int_0^tV^\omega_n-\langle V\rangle_nds\|_{L^\infty_{t, x}([0, T]\times\R^3}\|f_n\|_{S(0, T)}\\
&+\|\int_0^t\Delta V^\omega_n-\Delta\langle V\rangle_nds\|_{L^\infty_{t, x}([0, T]\times\R^3}\|f_n\|_{L^{\tilde q'}_tL^{\tilde r'}_x}\\
&+\|\int_0^t\nabla V^\omega_n-\nabla\langle V\rangle_nds\|_{L^\infty_{t, x}([0, T]\times\R^3}\|\nabla f_n\|_{L^{\tilde q'}_tL^{\tilde r'}_x}\\
&+\|\int_0^tV^\omega_n-\langle V\rangle_nds\|_{L^\infty_{t, x}([0, T]\times\R^3}\|(i\d_t+\Delta)f_n\|_{L^{\tilde q'}_tL^{\tilde r'}_x}.
\end{align*}
We notice that thanks to the properties of $V^\omega_n, \langle V\rangle_n$, we can apply Lemma \ref{lemma:unif_conv} to say that
\begin{equation*}
\|\int_0^tV^\omega_n-\langle V\rangle_nds\|_{L^\infty_{t, x}([0, T]\times\R^3)}\to0,\qquad\textrm{as}\;|\omega|\to\infty.
\end{equation*}
Furthermore, the same Lemma applies to any derivatives of $V^\omega_n, \langle V\rangle_n$, 
\begin{equation*}
\|\int_0^tD^\alpha V^\omega_n-D^\alpha\langle V\rangle_nds\|_{L^\infty_{t, x}([0, T]\times\R^3)}\to0,\qquad\textrm{as}\;|\omega|\to\infty,
\end{equation*}
for any multi-index $\alpha\in\N^3$.
\newline
Consequently, for each fixed $n\in\N$, the term $I_4$ converges to zero in $S(0, T)$ as $|\omega|\to\infty$.
\newline
Now, let $\eps>0$ be arbitrarily small, thus we can choose $n^\ast\in\N$ big enough so that 
\begin{equation*}
I_1+I_2+I_3\leq\eps,\qquad\forall\;n\geq n^\ast.
\end{equation*}
Thus, for the integral in \eqref{eq:conv}, we have
\begin{equation*}
\|\int_0^te^{i(t-s)\Delta}\pare{(V^\omega-\langle V\rangle)f}(s)ds\|_{S(0, T)}\leq\eps+I_4,
\end{equation*}
and in the limit as $|\omega|\to\infty$, we obtain
\begin{equation*}
\lim_{|\omega|\to\infty}\|\int_0^te^{i(t-s)\Delta}\pare{(V^\omega-\langle V\rangle)f}(s)ds\|_{S(0, T)}\leq\eps,\quad\forall\;\eps>0.
\end{equation*}
Hence the limit is zero and the Lemma is proved.
\newline\qed\newline
We are now ready to prove Theorem \ref{thm:main}.
\newline
Let $u^\omega, u$ be solutions to \eqref{eq:nls_omega}, \eqref{eq:nls_aver}, respectively, with the same initial datum, and let us consider the equation for the difference $v:=u^\omega-u$, which reads
\begin{equation}\label{eq:nls_diff}
\left\{\begin{array}{l}
i\d_tv=-\Delta v+(V^\omega-\langle V\rangle)u+V^\omega v
+C_1[(|\cdot|\ast|u^\omega|^2)u^\omega-(|\cdot|\ast|u|^2)u]
-a[|u^\omega|^\sigma u^\omega-|u|^\sigma u]\\
v(0)=0.
\end{array}\right.
\end{equation}
By the Duhamel's formula and by applying Strichartz estimates to each term we obtain
By applying Strichartz estimates we obtain
\begin{align*}
\|v\|_{S(0, T)}\lesssim&\|\int_0^te^{i(t-s)\Delta}\pare{(V^\omega-\langle V\rangle)u}(s)ds\|_{S(0, T)}\\
&+\|V^\omega_1v\|_{L^2_tL^{6/5}_x([0, T]\times\R^3)}
+\|V^\omega_2v\|_{L^1_tL^2_x([0, T]\times\R^3)}\\
&+\|(|\cdot|^{-1}\ast|u^\omega|^2)u^\omega-(|\cdot|^{-1}\ast|u|^2)u\|_{L^1_tL^2_x([0, T]\times\R^3}\\
&+\||u^\omega|^\sigma u^\omega-|u|^\omega u\|_{L^{\frac{4(\sigma+2)}{\sigma+8}}_tL^{\frac{\sigma+2}{\sigma+1}}_x}.
\end{align*}
Now we use inequalities \eqref{eq:due}, \eqref{eq:tre}, \eqref{eq:iii}, \eqref{eq:iv}, to get
\begin{align*}
\|v\|_{S(0, T)}\lesssim&\|\int_0^te^{i(t-s)\Delta}\pare{(V^\omega-\langle V\rangle)u}(s)ds\|_{S(0, T)}\\
&+T^{\frac{1-3\eps}{2}}\|V^\omega_1\|_{L^\infty_tL^{\frac{3}{1+3\eps}}_x}\|v\|_{S(0, T)}\\
&+T^{\frac{1+3\eps}{2}}\|V^\omega_2\|_{L^\infty_tL^{\frac{3}{1-3\eps}}_x}\|v\|_{S(0, T)}\\
&T^{\frac{1}{2}}\pare{\|u^\omega\|_{S(0, T)}^2+\|u\|_{S(0, T)}^2}\|v\|_{S(0, T)}\\
&T^{\frac{2(\sigma+2)}{4-\sigma}}\pare{\|u^\omega\|_{S(0, T)}^\sigma+\|u\|_{S(0, T)}^\sigma}\|v\|_{S(0, T)}\\
\lesssim&\|\int_0^te^{i(t-s)\Delta}\pare{(V^\omega-\langle V\rangle)u}(s)ds\|_{S(0, T)}
+T^\alpha C(\|u^\omega\|_{S(0, T)}, \|u\|_{S(0, T)})\|v\|_{S(0, T)}.
\end{align*}
Let us now recall that, by \eqref{eq:unif_bounds}, the norms $\|u^\omega\|_{S(0, T)}$ are uniformly bounded, indipendently on $\omega$. Consequently the constant $C$ in the inequality above depends only on $T$ and $\|u_0\|_{L^2(\R^3)}$. Thus we can use Gronwall's inequality to obtain
\begin{equation*}
\|v\|_{S(0, T)}\lesssim e^{cT}\|\int_0^te^{i(t-s)\Delta}\pare{(V^\omega-\langle V\rangle)u}(s)ds\|_{S(0, T)}.
\end{equation*}
By applying Lemma \ref{lemma:main} we finally prove that
\begin{equation*}
\|u^\omega-u\|_{S(0, T)}\to0,\qquad\textrm{as}\;|\omega|\to\infty,
\end{equation*}
for all finite times $0<T<\infty$. Hence the Theorem is proved.
\newline\qed

{\bf Acknowledgement:} The first and second author would like to thank the Department of Applied Mathematics and Theoretical Physics at the University of Cambridge for its hospitality and support during the preparation of this work.

\end{document}